\documentclass[11pt]{article}%

\newtheorem{theorem}{Theorem}
\newtheorem{lemma}{Lemma}
\newtheorem{proposition}{Proposition}
\newtheorem{corollary}{Corollary}

\newtheorem{assumption}{Assumption}

\newcommand{\lab}[1]{\label{#1}}
\newcommand{\labs}[1]{\label{#1}}
\newcommand{\labe}[1]{\label{#1}}

\newcommand{\dating}[3]{\date{{\Large {\bf #1}}\vspace{5mm}\\{\Large#3}}\label{#2}}

\RequirePackage[colorlinks,citecolor=blue,urlcolor=blue]{hyperref}

\usepackage{amsmath,amssymb,mathrsfs,eurosym}

\def\Bbb{\mathbb}
\def\scr{\mathscr}

\newcommand{\Section}[1]{\section{#1}\setcounter{equation}{0}}

\def\argmin{\mathop{\rm argmin}}

\newcommand{\etc}[1]{#1_1,\ldots,#1_n}

\newcommand{\cqfd}{\qquad\framebox[2.7mm]{\rule{0mm}{.7mm}}}
\newcommand{\bm}[1]{\mathbf{#1}}

\newcommand{\st}{\strut}

\begin{document}
\title{\mbox{}\vspace{-10mm}\\{\Huge{\bf About the non-asymptotic\vspace{2mm}\\behaviour of Bayes estimators\vspace{5mm}}} }
\author{
{\Large Lucien Birg\'{e}\vspace{3mm}}\\{\Large Sorbonne Universit\'es -- UPMC Universit\'e Paris 06\vspace{1mm}}\\
{\Large U.M.R.\  7599, L.P.M.A.\vspace{1mm}}\\
{\Large F-75005, Paris, France}
 }
\dating{\mbox{}}{May 2013}{September 2014}     

\maketitle

\begin{abstract}
This paper investigates the {\em nonasymptotic} properties of Bayes procedures for estimating an unknown distribution from $n$ i.i.d.\ observations. We assume that the prior is supported by a model $(\scr{S},h)$ (where $h$ denotes the Hellinger distance) with suitable metric properties involving the number of small balls that are needed to cover larger ones. We also require that the prior put enough probability on small balls. 

We consider two different situations. The simplest case is the one of a parametric model containing the target density for which we show that the posterior concentrates around the true distribution at rate $1/\sqrt{n}$. In the general situation, we relax the parametric assumption and take into account a possible mispecification of the model. Provided that the Kullback-Leibler Information between the true distribution and $\scr{S}$ is finite, we establish risk bounds for the Bayes estimators.
\end{abstract}

\Section{Introduction\labs{I}}
The purpose of this paper is to derive in a simple way some non-asymptotic results about posterior distributions and Bayes estimators from a frequentist viewpoint, therefore offering a complementary point of view to the classical results by Ghosal, Ghosh and van der Vaart (2000) --- see also the related papers: Ghosal, Lember and van der Vaart (2003) and van der Vaart (2003) ---. It can also be considered as a new and extended presentation of Le Cam (1973 and 1982). In any case, it has been strongly influenced by these three papers.

We shall work here within the following framework: we have at disposal a sample 
$\bm{X}=(\etc{X})$ of size $n$, the $X_i$ being measurable mappings from 
$(\Omega,{\cal A})$ to $({\scr X},{\cal X})$ with a common unknown distribution $P$. This distribution is an element of the metric space $({\scr P},h)$ of all probability measures on $({\scr X},{\cal X})$ endowed with the Hellinger distance $h$ given by
\[
h^2(R,T)=\frac{1}{2}\int\left(\sqrt{\frac{dR}{d\lambda}}-\sqrt{\frac{dT}{d\lambda}}\right)^2
d\lambda,
\]
where $\lambda$ is an arbitrary positive measure which dominates both $R$ and $T$. We then introduce a model for $P$, i.e.\ a dominated family ${\scr S}=\{P_t\,|\,t\in S\}\subset\scr{P}$ of probabilities on ${\scr X}$ with densities $f_t=dP_t/d\mu$ with respect to some reference measure 
$\mu$ on ${\scr X}$. We assume that the mapping $t\mapsto P_t$ is one to one which allows us to systematically identify $S$ and ${\scr S}$, thus considering $S$ as a metric space with distance 
$h$ --- $h(t,u)=h(P_t,P_u)$ --- and the corresponding Borel $\sigma$-algebra. We then introduce a prior distribution $\nu$ on $S$, turning the parameter $t$ into a random variable $\bm{t}$. The prior $\nu$ and the sample $\bm{X}$ give rise to a posterior distribution $\overline{\nu}=\overline{\nu}(\cdot|\bm{X})$ and, given a loss function $w\circ h$ on $S\times S$, to a corresponding Bayes estimator $\widetilde{s}$ defined by
\begin{equation}
\widetilde{s}(\cdot|\bm{X})=\argmin_{u\in S}\Bbb{E}\left[w\left(h(u,\bm{t})\right)\!|\bm{X}\right]=\int_{S}w\left(h(u,\bm{t})\right)d\overline{\nu}({\bm t|{\bm X}}),
\labe{Eq-Bayesest}
\end{equation}
where $\argmin$ refers to any minimizer in case it is not unique. In the sequel we shall write $\Bbb{E}_s[f(\bm X)]$ to indicate that the $X_i$ are i.i.d.\ with distribution $P_s$ and $\Bbb{P}_s$ for the corresponding probability on $\Omega$ that gives $\bm{X}$ the distribution $P_s^{\otimes n}$.

Our purpose here will be twofold. When $P=P_s$ truly belongs to ${\scr S}$ and the metric structure of 
$(\scr{S},h)$ is similar to that of a compact subset of some Euclidean space, we shall study the concentration rate of the posterior distribution $\overline{\nu}(\cdot|\bm{X})$ of $\bm{t}$ around $P_s$. When $P$ does not belong to ${\scr S}$ or the metric structure of $\scr{S}$ does not follow the previous requirements, we shall study the performance of the Bayes estimator(s) $P_{\widetilde{s}}$ of $P$ defined via the loss function $w\circ h$ for suitable functions $w$.
The main feature of our approach is its non-asymptotic viewpoint, explicit deviation bounds being provided for fixed $n$.

\paragraph{Some notations}
To begin with, let us fix some notations to be used throughout the paper.
In the metric space $(S,h)$, we denote by ${\cal B}(t,r)$ the closed Hellinger ball with center 
$t\in S$ and radius $r$ while the ball with center $y$ and radius $r$ in the Euclidean space 
$\Bbb{R}^d$ will be denoted ${\cal B}_d(y,r)$. The set of positive integers 
$\Bbb{N}\setminus\{0\}$ is denoted by $\Bbb{N}^*$, the cardinality of the set $N$ by $|N|$ and 
we write $a\vee b$ for $\max\{a,b\}$. The distance between two sets $A$ and $B$ is 
$h(A,B)=\inf_{t\in A,\,u\in B}h(t,u)$ and if $\bm{x}=(\etc{x})\in{\scr X}^n$ we write $f_t(\bm{x})$ instead of $\prod_{i=1}^nf_t(x_i)$. 

For any measurable subset $B$ of $S$ such that $\nu(B)>0$, we define the density $g_B$ with respect to $\mu^{\otimes n}$ and the probability $P_B$ on ${\scr X}^n$ by
\begin{equation}
g_B(\bm{x})=\frac{1}{\nu(B)}\int_Bf_{t}(\bm{x})\,d\nu(t)\qquad\mbox{and}\qquad
P_B=g_B\cdot\mu^{\otimes n}.
\labe{Eq-gB}
\end{equation}
We denote by $\Bbb{P}_B$ the probability on $\Omega$ that gives $\bm{X}$ the distribution $P_B$.

\Section{A toy example\labs{T}}
Let us first consider, in order to motivate our approach, the very particular situation of a finite or countable parameter set $S$ containing the true density $s$ to estimate. Besides, we shall assume that $\nu(t)>0$ for all $t\in S$. The posterior probability is given in this case by
\[
\overline{\nu}(B|{\bm X})=\frac{\sum_{t\in B}\nu(t)f_t(\bm{X})}{\sum_{t\in S}\nu(t)f_t(\bm{X})}=\left(1+{\sum_{t\in B^{c}}\nu(t)f_t(\bm{X})\over \sum_{t\in B}\nu(t)f_t(\bm{X})}\right)^{-1}\ge 1-\frac{\sum_{t\in B^c}\nu(t)f_t(\bm{X})}{\sum_{t\in B}\nu(t)f_t(\bm{X})}
\]
for all $B\subset S$ and it follows that
\begin{equation}
\overline{\nu}(B|{\bm X})\ge 1-\frac{\sum_{t\in B^c}\nu(t)f_t(\bm{X})}{\nu(s)f_s(\bm{X})}
\quad\mbox{for all }B\ni s.
\labe{Eq-aux9}
\end{equation}
In order to evaluate the concentration of the posterior distribution $\overline{\nu}(.|{\bm X})$ around $s$, we focus on those sets $B_k$ which are Hellinger balls centered at $s$ with radius $k/\sqrt{n}$, $k\in\Bbb{N}^*$. Bounding $\overline{\nu}(B_k|{\bm X})$ from below requires to bound from above ratios of the form $f_t(\bm{X})/f_s(\bm{X})$ when the $X_i$ are distributed according to $P_s$, which implies that $f_s(\bm{X})>0$ a.s. This control derives from Lemma~7 in Birg\'e (2006) which implies the
following inequality:
%
\begin{lemma}\lab{L-expineq}
Given $n$ i.i.d.\ random variables $\etc{X}$ with distribution $P$ and another distribution $Q$, 
then $\log\left((dQ/dP)(X_i)\st\right)\in[-\infty,+\infty)$ a.s.\ (with the convention $\log0=-\infty$) and, for all $y\in\Bbb{R}$,
\[
\Bbb{P}\left[\sum_{i=1}^n\log\left(\frac{dQ}{dP}(X_i)\right)\ge y\right]\le\exp\left[-\frac{y}{2}\right]
\rho^n(P,Q)\quad\mbox{with}\quad\rho(P,Q)=\int\sqrt{\frac{dP}{d\lambda}\frac{dQ}{d\lambda}}\,
d\lambda.
\]
\end{lemma}
We recall here that $\rho(P,Q)$ is called the {\em Hellinger affinity} between $P$ and $Q$, the definition being independent of the choice of the dominating measure $\lambda$, and that it satisfies
\begin{equation}
\rho(P,Q)=1-h^2(P,Q)\quad\mbox{and}\quad\rho\left(P^{\otimes n},Q^{\otimes n} \right)=\rho^n(P,Q)\le\exp\left[-nh^2(P,Q)\right],
\labe{Eq-rhon}
\end{equation}
hence
\begin{equation}
h^2\left(P^{\otimes n},Q^{\otimes n} \right)=1-\rho^n(P,Q)=1-\left(1-h^2(P,Q)\right)^n\le nh^2(P,Q).
\labe{Eq-hn}
\end{equation}
We therefore derive from Lemma~\ref{L-expineq} and (\ref{Eq-rhon}) that, for $\delta>0$,
\[
\Bbb{P}_{s}\left[f_t(\bm{X})\ge\delta\nu(s)f_s(\bm{X})\right]\le[\delta\nu(s)]^{-1/2}
\rho^n(P_s,P_t)\le[\delta\nu(s)]^{-1/2}\exp\left[-nh^2(P_s,P_t)\right].
\]
Setting, for $k\in\Bbb{N}^*$,
\[
\Gamma_k=\left\{\bm{x}\,\left|\,\sup_{t\in B_k^c}f_t(\bm{x})\ge
\delta\nu(s)f_s(\bm{x})\right.\right\}
\]
and denoting by $N_l$, for $l\in\Bbb{N}^*$, the cardinality of the set
\[
B_{l+1}\setminus B_l=\left\{t\in S\,\left|\ {l\over \sqrt{n}}< h(P_{s},P_{t})\le {l+1\over \sqrt{n}}\right.\right\},
\]
we conclude that
\[
\sqrt{\delta\nu(s)}\,\Bbb{P}_s[\bm{X}\in\Gamma_k]\le\sum_{t\in B_k^c}
\exp\left[-nh^2(P_s,P_t)\right]\le\sum_{l\ge k}\sum_{t\in B_{l+1}\setminus B_l}e^{-l^2}
=\sum_{l\ge k}N_le^{-l^2}.
\]
Consequently, if the series $\sum_{l\ge 1}N_le^{-l^2}$ converges, $\Bbb{P}_s[\bm{X}\in\Gamma_k]$ can be made smaller than any $\epsilon>0$, provided that $k=k(\epsilon,\delta,\nu(s))$ is large enough (note that the smaller $\nu(s)$ or $\delta$, the larger $k(\epsilon,\delta,\nu(s))$). For such a $k$, 
$\Bbb{P}_s[\bm{X}\in\Gamma_k^c]\ge1-\epsilon$ and (\ref{Eq-aux9}) implies that $\overline{\nu}(B_k|{\bm X})\ge1-\delta\nu(B^c_k)\ge1-\delta$ provided that ${\bm X}\in\Gamma_k^c$. These conclusions can be summarized in the following way:
%
\begin{proposition}\label{P-toystory2}
Consider the previous toy example where $P=P_s$ and assume that the series $\sum_{l\ge 1}
N_le^{-l^2}$ converges. For all $s\in S$ and $\epsilon,\delta\in(0,1)$, with probability at least $1-\varepsilon$ with respect to $P_s^{\otimes n}$, 
\[
\overline \nu\left[\left.B\left(s,k/\sqrt{n}\right)\right|{\bm X}\right]\ge 1-\delta,
\]
where 
\[
k=k(\epsilon,\delta,\nu(s))=\inf\left\{j\ge 1\,\left|\,\sum_{l\ge j}N_le^{-l^2}\le\epsilon\sqrt{\delta\nu(s)}\right.\right\}
\]
is independent of $n$ and nonincreasing with the positive parameters $\epsilon$, $\delta$ and $\nu(s)$.
\end{proposition}
Of course, the situation of interest is more general than that of our toy example and one would like to extend this result not only to more general sets $S$ but also to cases where the true distribution $P$ does not belong to ${\scr S}$. The main idea that allows to handle these extensions is to replace $P$ by some $P_s$ that belongs to ${\scr S}$, $s$ by a small Hellinger ball around it and the countable parameter set of our toy example by a countable covering of $S$ by small balls. 

\Section{Concentration of the posterior distribution\labs{A}}

In the sequel, we systematically fix an arbitrary point in $S$, denoted by $s$, which may or may not satisfy $P=P_s$. One should think of $P_s$ as an approximation for $P$ although it can be any point in $S$. It will remained fixed throughout the paper and play a central role for the statements of our assumptions and results since all of them will be related to it and involve balls centered at $s$. Therefore, given $k\in\Bbb{N}$, we shall write $\scr{B}(k)$ instead of 
${\cal B}(s,2^{-k})$ for simplicity. 

We also want to recall some well-known facts about metric entropies and related notions that can, for instance, be found in Kolmogorov and Tikhomirov (1961) or Lorentz (1966). Let $(M,d)$ be a metric space. A subset $N$ of $M$ is called $x$-separated (with $x>0$) if any two distinct points in $N$ are at a distance larger than $x$. Given a subset $A$ of $M$ we denote by $N_A(x)$ the smallest number of closed balls of radius $x>0$ in $M$ that are needed to cover $A$.
%
\begin{lemma}\lab{L-entro}
Let $A$ be a subset of a metric space $(M,d)$ and $N$ be a maximal $x$-separated subset $N$ of $A$.
Then $N_A(x)\le|N|\le N_A(x/2)$.
\end{lemma}
%
\subsection{Our assumptions\labs{A2}}
In our toy example, the concentration of the posterior distribution around $P_s$ depends on the amount of mass that the prior puts on the true distribution $P_{s}=P$ and on the number of points that are contained in Hellinger balls centered on it. In the general situation, the point $s$ has to be replaced by a small ball around $s$ and our assumptions have to be modified accordingly. The first one deals with the number of small balls that are needed to cover larger ones as an alternative to the behaviour of the numbers $N_l$ of our toy example.
%
\begin{assumption}\lab{A-Dim}
There exists a nonincreasing function $D$ from $(0,1/4]$ to $[1,+\infty)$ such that, for all 
$x\in(0,1/4]$, any closed ball ${\cal B}(\cdot,4x)$ with radius $4x$ in $S$ and any $x$-separated subset $N$ of ${\cal B}(\cdot,4x)$, $|N|\le\exp[D(x)]$. 
\end{assumption}
Such an assumption holds for Euclidean balls in $\Bbb{R}^d$ with a function $D$ bounded by $d\log9$ as can be checked in the following way using a volume argument: all Euclidean balls of radius $x/2$ with centers in $N$ are disjoint and included in a ball of radius $9x/2$ so that their number is bounded by $9^d$. This result can be used to check Assumption~\ref{A-Dim} when $S\subset\Bbb{R}^d$, which means that $\scr{S}$ is a parametric model, and there exists a polynomial relationship in $S$ between the Hellinger and Euclidean distances of the form
\begin{equation}
a\|t-u\|^\alpha\le h(t,u)\le A\|t-u\|^\alpha,\quad\mbox{for some }A>a>0,\;\alpha>0
\mbox{ and all }t,u\in S.
\labe{Eq-heleuc}
\end{equation}
More generally, checking Assumption~\ref{A-Dim} amounts to making entropy computations as shown by Lemma~\ref{L-entro}. Indeed, since
\[
h(t,u)=\frac{1}{\sqrt{2}}\left\|\sqrt{f_t}-\sqrt{f_u}\right\|_2\quad\mbox{with}\quad
\left\|f-g\right\|_2^2=\int(f-g)^2d\mu,
\]
one can use the mapping $t\mapsto\psi(t)=\sqrt{f_t}$ from $S$ into $\Bbb{L}_2(\mu)$ to derive entropy bounds for $S$ from entropy bounds for $\psi(S)$ in $\Bbb{L}_2(\mu)$. For suitable models $S$, the images $\psi(S)$ are compact subsets of function spaces (like H\"older and Besov spaces) for which these entropy bounds are well-known, typically resulting in functions of the type $D(x)=cx^{-\delta}$, $c,\delta>0$.

Assumption~\ref{A-Dim} also provides a bound for $N_S\left(4^{-j}\right)$ for all $j\ge1$. Indeed, to cover $S$ with balls of radius $4^{-k}$ with $k\ge1$ it suffices to cover $S$ with balls of radius $4^{-k+1}$ and then cover each of these balls of radius $4\times4^{-k}$ with smaller ones of radius $4^{-k}$. Lemma~\ref{L-entro} and Assumption~\ref{A-Dim} imply that $N_{{\cal B}(\cdot,4x)}(x)\le\exp[D(x)]$ for any ball 
${\cal B}(\cdot,4x)$ with radius $4x\le1$, which leads to the recursive formula
\begin{equation}
N_{S}\left(4^{-k}\right)\le N_{S}\left(4^{-k+1}\right)\exp\left[D\left(4^{-k}\right)\right]
\quad\mbox{for all }k\ge1.
\labe{Eq-rec}
\end{equation}
Finally since $S={\cal B}(t,1)={\cal B}(t,4\times 4^{-1})$ for all $t\in S$ and $D(\cdot)$ is nonincreasing, it follows from (\ref{Eq-rec}) that
\[
N_{S}\left(4^{-j}\right)\le\prod_{k=1}^{j}\exp\left[D\left(4^{-k}\right)\right]\le
\exp\left[jD\left(4^{-j}\right)\right]\quad\mbox{for all }j\ge 1.
\]

Our second assumption is about the behaviour of the prior $\nu$ in vicinities of $s$.
\begin{assumption}\lab{A-Prior}
Let $\scr{B}(j)$ denote the closed Hellinger ball with center $s$ and radius $2^{-j}$ for 
$j\in\Bbb{N}$. There exists a nondecreasing function $\beta$ from $\Bbb{N}$ to $[1,+\infty)$ such that the following property holds for the prior $\nu$ and some $\gamma\in[1,4]$: 
\begin{equation}
\nu\!\left(\st{\scr B}(j-k)\right)\le\exp\left[\gamma^{k}\beta(j)\right]\nu\!\left(\st{\scr B}(j)\right)
\quad\mbox{for all }j\ge3\quad\mbox{and}\quad 3\le k\le j.
\labe{Eq-AS2}
\end{equation}
\end{assumption}
Here are some comments in order. First of all, if the assumption holds for some $\gamma<4$ it also holds with $\gamma=4$. Then, on the one hand, (\ref{Eq-AS2}) with $k=j$ requires that 
$\nu\!\left(\st{\scr B}(j)\right)\ge\exp\left[-\gamma^{j}\beta(j)\right]$ since $\nu({\scr B}(0))=\nu(S)=1$.
On the other hand, since $k\ge3$ and $\gamma\ge1$, the right-hand side of (\ref{Eq-AS2}) is not smaller than $\exp\left[\gamma^3\beta(j)\right]\nu\!\left(\st{\scr B}(j)\right)$ and the assumption is satisfied as soon as $\nu\!\left(\st{\scr B}(j)\right)\ge\exp\left[-\gamma^3\beta(j)\right]$.

This second assumption will in particular be typically satisfied by many parametric models 
$\scr{S}$ when $S$ is a bounded subset of $\Bbb{R}^d$ and $\nu$ the uniform distribution on $S$. Indeed, the Lebesgue measure $\lambda$ on $\Bbb{R}^d$ satisfies an analogue of (\ref{Eq-AS2}), namely
\[
\lambda\!\left(\st{\cal B}_d(\cdot,2^{k-j})\right)=\exp\left[kd\log 2\right]
\lambda\!\left(\st{\cal B}_d(\cdot,2^{-j})\right)
\le\exp\left[(3/2)^{k}(d\log 2)\right]\lambda\!\left(\st{\cal B}_d(\cdot,2^{-j})\right),
\]
for $3\le k\le j$, where ${\cal B}_d(\cdot,r)$ denotes an arbitrary  ball of radius $r$ in 
$\Bbb{R}^d$. As a consequence, this property still holds when $\lambda$ is replaced by the uniform distribution $\nu$ on the bounded set $S$ and $d\log2$ by another suitable constant, provided that the shape of $S$ is regular enough to imply that the ratio between the volume of ${\cal B}_d(t,r)$ and the volume of ${\cal B}_d(t,r)\cap S$ is bounded from above uniformly for $t\in S$. If, moreover, (\ref{Eq-heleuc}) holds, Assumption~\ref{A-Prior} will be satisfied for the uniform prior on $S$ with a constant function $\beta$.

\subsection{A basic concentration result\labs{A3}}
The following important auxiliary result to be proved in Section~\ref{P} below is not exactly 
what we would like to get since it involves the probability $P_{{\scr B}(J_1)}$ (as defined by (\ref{Eq-gB}) with ${\scr B}(J_1)$ the Hellinger ball with center $s$ and radius $2^{-J_1}$) instead of the true distribution $P^{\otimes n}$ of $\bm{X}$. It will nevertheless be the main 
tool for our treatment of the cases of interest.
%
\begin{theorem}\lab{T-Bayes}
Let $\overline{\nu}(\cdot|\bm{X})$ be the posterior distribution of $\bm{t}$, $s$ be a given element of $S$ and $\scr{B}(j)={\cal B}(s,2^{-j})$ for $j\in\Bbb{N}$. Let 
Assumptions~\ref{A-Dim} and \ref{A-Prior} hold and 
\begin{equation}
n\ge4^4\left(\left[\gamma^{3}\beta(3)/3\right]\bigvee D(1/4)\right).
\labe{Eq-n}
\end{equation}
Then one can find integers $J_1$ and $J$ satisfying
\begin{equation}
1\le J\le J_1-2\qquad\mbox{and}\qquad
4^{-(J+3)}n\ge\left[\gamma^{J_1-J+1}\beta(J_1)/3\right]\bigvee D\!\left(2^{-(J+1)}\right).
\labe{Eq-Required}
\end{equation}
For any pair $(J_1,J)$ for which (\ref{Eq-Required}) holds, there exists some $\Gamma_J\subset\scr{X}^n$ with
\begin{equation}
P_{{\scr B}(J_1)}[\Gamma_J]=\int_{\Gamma_J}g_{{\scr B}(J_1)}(\bm{x})\,d\mu^{\otimes n}(\bm{x})\le
1.05\exp\left[-4^{-(J+3)}n\right]
\labe{Eq-T1a}
\end{equation}
and such that, if $\bm{X}\not\in\Gamma_J$, the posterior distribution satisfies
\begin{equation}
\overline{\nu}\left[\left.{\scr B}\left(j\right)\right|\bm{X}\right]\ge1-1.05\exp\left[-4^{-(j+3)}n\right]
\quad\mbox{for all }\:j\in\Bbb{N}, \;j\le J.
\labe{Eq-T1b}
\end{equation}
\end{theorem}
Since $D$ is nonincreasing, bounded from below by one and $\gamma\le4$, (\ref{Eq-Required}) also implies that $4^{-(J+3)}n\ge1$ and
\begin{equation}
4^{-(J'+3)}n\ge\left[\gamma^{J_1-J'+1}\beta(J_1)/3\right]\bigvee D\!\left(2^{-(J'+1)}\right)
\quad\mbox{for }J'\in\Bbb{N},\;\;1\le J'\le J.
\labe{Eq-Req}
\end{equation}
This means that (\ref{Eq-T1a}) and (\ref{Eq-T1b}) also hold with $J'<J$ replacing $J$ 
everywhere. Choosing such a $J'$ allows to find a set $\Gamma_{J'}$ for which 
$P_{{\scr B}(J_1)}[\Gamma_{J'}]\le1.05\exp\left[-4^{-(J'+3)}n\right]$, which improves the 
bound in (\ref{Eq-T1a}).

\Section{The case of a ``true" model\labs{B}}
Throughout this section, we shall assume that the model $S$ is ``almost true" which means 
that our observations are distributed according to some distribution $P$ which is very close to
$P_s\in\scr{S}$.
%
\begin{theorem}\lab{T-truemodel}
Let Assumptions~\ref{A-Dim} and \ref{A-Prior} hold with
\[
\gamma<4,\qquad\overline{\beta}=\sup_{j\ge3}\beta(j)<+\infty\qquad\mbox{and}
\qquad\overline{D}=\sup_{0<x\le1/4}D(x)<+\infty.
\]
Let $c\in(0,1/2]$ and assume that
\begin{equation}
n\ge4^{4}\left(\left[\left(16c^{-2}n\right)^{(\log\gamma)/(\log4)}\overline{\beta}/3\right]\bigvee\overline{D}\bigvee\log(8/c)\right)
\labe{Eq-J0}
\end{equation}
and that the true distribution $P$ of the $X_i$ satisfies $h(P_s,P)=\kappa/\sqrt{2n}$ with 
$0\le\kappa<1/2$. For $x>0$, let $J_2(x)$ be the smallest positive integer $j$ such that
\begin{equation}
4^{j-4}x^2\ge\left(\gamma^{j+2}\overline{\beta}/3\right)\bigvee\overline{D}\bigvee\log(8/x).
\labe{Eq-J4}
\end{equation}
Then, with probability at least $1-(\kappa+c)$ with respect to $P^{\otimes n}$,
\begin{equation}
\overline{\nu}\left[\left.{\cal B}\!\left(s,2^{k}c/\sqrt{n}\right)\right|\bm{X}\right]\ge
1-1.05\exp\left[-4^{k-4}c^2\right]\quad\mbox{for all integers }k\ge J_2(c).
\labe{Eq-J1}
\end{equation}
%
\end{theorem}
%
{\em Proof:} First note that the function $J_2$ is well-defined and that (\ref{Eq-J0}) holds for 
$n$ large enough since $\gamma<4$. Let $c_n\in[c,2c)$ be such that 
$\log\left(2\sqrt{n}/c_n\right)/\log2\in\Bbb{N}$ and 
$J_1=\left[\log\left(\sqrt{n}/c_n\right)/\log2\right]+2$; then $c_n=2^{-J_1+2}\sqrt{n}$ and 
$J_1\ge7$ since $c_n<1$ and $n\ge4^4$. Moreover, since $n=4^{J_1-2}c_n^2$ and (\ref{Eq-J0}) remains true with $c_n$ replacing $c$, it implies that
\[
4^{J_1-6}c_n^2\ge\left[\left(4^{J_1}\right)^{(\log\gamma)/(\log4)}\overline{\beta}/3\right]
\bigvee\overline{D}\bigvee\log(8/c_n),
\]
which means that (\ref{Eq-J4}) is satisfied for $x=c_n$ and $j=J_1-2$. Therefore $J_0=J_2(c_n)$ satisfies
\[
1\le J_0\le J_1-2\qquad\mbox{and}\qquad4^{J_0-4}c_n^2\ge\left(\gamma^{J_0+2}\overline{\beta}/3\right)\bigvee\overline{D}\bigvee\log(8/c_n).
\]
or, equivalently,
\begin{equation}
4^{J_0-2-J_1}n\ge\left(\gamma^{J_0+2}\overline{\beta}/3\right)\bigvee\overline{D}
\bigvee\log(8/c_n).
\labe{Eq-J6}
\end{equation}
This means that (\ref{Eq-Required}) holds with $J=J_1-J_0-1$ and Theorem~\ref{T-Bayes} applies.
We first derive from (\ref{Eq-T1a}) and (\ref{Eq-J6}) that
\[
P_{{\scr B}(J_1)}[\Gamma_J]\le1.05\exp\left[-4^{-(J_1-J_0+2)}n\right]=
1.05\exp\left[-4^{J_0-4}c_n^2\right]\le1.05c_n/8.
\]
Since, for $t\in{\scr B}(J_1)$, $h(t,s)\le2^{-J_1}$, it follows that $h(P,P_t)\le 
c_n/\!\left(4\sqrt{n}\right)+\kappa/\sqrt{2n}$ hence $h^2\left(P_{t}^{\otimes n},P^{\otimes n}\right)\le\left(\kappa/\sqrt{2}+c_n/4\right)^2$ by (\ref{Eq-hn}). Then Lemma~\ref{L-bayconv} below, to be proved in Section~\ref{P}, shows that $\rho\left(P_{{\scr B}(J_1)},
P^{\otimes n}\right)\ge1-\left(\kappa/\sqrt{2}+c_n/4\right)^2$. Hence, by a classical result of Le Cam (1973) relating the variation distance to the Hellinger affinity,
\begin{eqnarray*}
\sup_A\left[P_{{\scr B}(J_1)}(A)-P^{\otimes n}(A)\right]^2&\le&1-
\rho^2\left(P_{{\scr B}(J_1)},P^{\otimes n}\right)\;\;\le\;\;
1-\left[1-\left(\kappa/\sqrt{2}+c_n/4\right)^2\right]^2\\&<&2\left(\kappa/\sqrt{2}+c_n/4\right)^2
\end{eqnarray*}
and finally,
\[
\Bbb{P}[\bm{X}\in\Gamma_J)]=P^{\otimes n}\left(\Gamma_J\right)<1.05(c_n/8)+
\sqrt{2}\left(\kappa/\sqrt{2}+c_n/4\right)<\kappa+(c_n/2)<\kappa+c.
\]
Now observe that (\ref{Eq-T1b}) can be written
\[
\overline{\nu}\left[\left.{\cal B}\!\left(s,2^{J_1-2-j}c_n/\sqrt{n}\right)\right|\bm{X}\right]\ge
1-1.05\exp\left[-4^{J_1-j-5}c_n^2\right]\quad\mbox{for all }\,j\in\Bbb{N}, \;j\le J.
\]
Setting $k=J_1-j-1$ leads to 
\begin{equation}
\overline{\nu}\left[\left.{\cal B}\!\left(s,2^{k-1}c_n/\sqrt{n}\right)\right|\bm{X}\right]\ge
1-1.05\exp\left[-4^{k-4}c_n^2\right]\quad\mbox{for }J_0\le k\le J_1-1,
\labe{Eq-J5}
\end{equation}
while, for larger values of $k$, $2^{k-1}c_n/\sqrt{n}>1$ and (\ref{Eq-J5}) still holds. Finally, (\ref{Eq-J1}) follows from (\ref{Eq-J5}) and the monotonicty of $J_2$.\cqfd
%
\begin{lemma}\lab{L-bayconv}
Let $\lambda$ be a probability on $S$ such that for some $P\in{\scr P}$,
$\lambda(\{t\,|\,h(P,P_t)\le r\})=1$ and let $P_{\lambda}$ be the distribution with density
$f_\lambda(x)=\int_S f_t(x)\,d\lambda(t)$ with respect to $\mu$. Then,
$\rho(P,P_\lambda)\ge 1-r^2$ or, equivalently, $h(P,P_\lambda)\le r$.
\end{lemma}
%

\paragraph{Interpretation}
As we already mentioned in Section~\ref{A2}, the case of bounded functions $\beta$ and $D$ typically occurs in a parametric situation when $S$ is a bounded subset of $\Bbb{R}^d$ for which (\ref{Eq-heleuc}) holds and $\nu$ is the uniform prior, although more exotic priors can also be considered since they only have to satisfy (\ref{Eq-AS2}). The choice of $c$, provided that $n$ is large enough to satisfy (\ref{Eq-J0}), allows to get a control of the probability of the event $\bm{X}\not\in\Gamma_J$, for which (\ref{Eq-J1}) holds, by $1-c-\kappa$. As to (\ref{Eq-J1}), it provides concentration properties of the posterior distribution for balls with center $s$ and radius at least $c2^{J_2(c)}/\sqrt{n}$. This confirms that, in the parametric case with $P=P_s$ and under a suitable assumption on the prior, the posterior concentrates around the true value $s$ at rate $n^{-1/2}$ (with respect to the Hellinger distance) as shown by Le Cam (1973) and Ibragimov and Has'minskii (1981). For parametric models with $S\in\Bbb{R}^d$ such that (\ref{Eq-heleuc}) holds, we can also derive concentration rates for the posterior with respect to the Euclidean distance.

\Section{Convergence of Bayes estimators\labs{C}}
If the requirements of Section~\ref{B} are not satisfied, we may still apply Theorem~\ref{T-Bayes} with $J_1=J_1(n)$ satisfying
\[
J_1(n)\ge3\qquad\mbox{and}\qquad4^{-J_1(n)-1}n\ge\left[\gamma^{3}\beta(J_1(n))/3\right]\bigvee D\!\left(2^{1-J_1(n)}\right)
\]
which implies by monotonicity that (\ref{Eq-Required}) holds for $0\le J\le J_1(n)-2$. It is always possible to find such an integer $J_1(n)$, at least for $n$ large enough. If $D\vee\beta$ is bounded we can set $2^{-J_1(n)}=cn^{-1/2}$ for a suitably large value of the constant $c$ as we did in the previous section. Otherwise $2^{-J_1(n)}\sqrt{n}$ goes to infinity with $n$, but, since $\beta$ and $D$ are nonincreasing functions, one can choose $2^{-J_1(n)}$ converging to zero when $n$ goes to infinity. It then follows from Theorem~\ref{T-Bayes} that, except on a set $\Gamma_J\subset\scr{X}^n$ with $P_{{\scr B}(J_1(n))}$-probability bounded by 
$1.05\exp\left[-4^{-J-3}n\right]$, the posterior distribution satisfies
\begin{equation}
\overline{\nu}\left[\left.{\scr B}\left(j\right)\right|\bm{X}\right]\ge1-1.05\exp\left[-4^{-(j+3)}n\right]
\quad\mbox{for all }\:j\in\Bbb{N}, \;j\le J.
\labe{Eq-aux1}
\end{equation}
Unfortunately, in such a case, we cannot get a control on $P^{\otimes n}(\Gamma_J)$ since we are unable to bound $\left|P_{{\scr B}(J_1(n))}(\Gamma_J)-P^{\otimes n}(\Gamma_J)\right|$ in a non-trivial way. Therefore, instead of studying the concentration properties of $\overline{\nu}$, we shall content ourselves to use a different approach which only provides concentration properties for Bayes estimators.

\subsection{Preliminary considerations\labs{C1}}
In order to define a Bayesian estimator according to (\ref{Eq-Bayesest}) we first have to fix our loss function $w\circ h$ and therefore to choose a function $w$. It should satisfy the following requirements for all $z\in(0,1/2]$ and suitable positive constants $\delta,a',B'$:
\begin{equation}
w(0)=0\quad\mbox{and}\quad x^\delta w(z)\le  w(xz)\le a'\exp\left[B'x^2\right]w(z)\quad\mbox{for }2\le x\le z^{-1}.
\labe{Eq-loss}
\end{equation}
Functions $w(z)=z^\delta$ with $\delta>0$ do satisfy (\ref{Eq-loss}) with $B'>0$ provided that $\delta\log x\le\log a'+B'x^2$ for $x\ge2$, which holds with $a'=\sup_{x\ge2}\exp\left[\delta\log x-B'x^2\right]$. But there are other less trivial cases as shown by the following proposition.
%
\begin{proposition}\lab{P-loss function}
For $0<z\le1/2$, the function $w(z)=\exp[\theta z^\delta]-1$ with $\theta>0$ and $0<\delta\le2$ satisfies (\ref{Eq-loss}) provided that $B'\ge\theta$ and $a'B'\ge1$. It holds in particular with $B'=\theta\vee1$ and $a'=1/B'$.
\end{proposition}
%
{\em Proof:}
On the one hand, 
\[
w(xz)=\exp\left[\theta x^\delta z^\delta\right]-1\ge x^\delta\left(\exp\left[\theta z^\delta\right]-1\right)=
x^\delta w(z),
\]
since $a(e^b-1)\le e^{ab}-1$ for $b\ge0$ and $a\ge1$. On the other hand, since 
$w(z)>\theta z^\delta$, it is enough to prove that $\exp[\theta(xz)^\delta]-1\le a'\theta z^\delta
\exp\left[B'x^2\right]$. Expending both sides in series we see that the inequality is satisfied provided that
$\left(\theta(xz)^\delta\right)^k\le a'\theta z^\delta\left(B'x^2\right)^k$ for all $k\ge1$ or, equivalently,
$1\le a'B'x^{(2-\delta)k}\left(B'\theta^{-1}\right)^{k-1}z^{-\delta(k-1)}$ which holds under our assumptions.\cqfd\vspace{2mm}

Since our prior $\nu$ is concentrated on the model $\scr{S}$ but we do not assume that 
the true distribution of the $X_i$ does belong to the model $\scr{S}$, we need a control of the approximation properties of $P$ by $\scr{S}$ since it is clear that there is no hope to get anything useful if we have chosen a model which does not approximate $P$ closely enough. But, unlike for universally robust estimators like the T-estimators of Birg\'e (2006) and the $\rho$-estimators of Baraud, Birg\'e and Sart (2014), we have to measure the distorsion between $\scr{S}$ and $P$ not in terms of the Hellinger distance, as would be the case for the previous estimators, but in terms of the Kullback-Leibler divergence, as is the case for the maximum likelihood estimator --- see for instance Massart (2007) ---.
Assumptions involving the Kullback-Leibler divergence also appear in other works on Bayesian estimators like Ghosal, Gosh and van der Vaart (2000). We recall that the Kullback-Leibler divergence between two probabilities $P,Q$ is given 
by 
\[
K(P,Q)=\int\log\left(dP\over dQ\right)dP\in[0,+\infty]\;\;\mbox{if }P\ll Q\quad\mbox{and}\quad
K(P,Q)=+\infty\;\;\mbox{otherwise}.
\]

\subsection{Application to Bayesian estimators\labs{C2}}
We are now in a position to prove the following theorem which describes the properties of Bayesian estimators.
%
\begin{theorem}\lab{T-Bayesest}
Let Assumptions~\ref{A-Dim} and \ref{A-Prior} hold and $n$ and $J_1(n)$ satisfy
\begin{equation}
J_1(n)\ge3\qquad\mbox{and}\qquad4^{-J_1(n)-1}n\ge\left[\gamma^{3}\beta(J_1(n))/3\right]\bigvee D\!\left(2^{1-J_1(n)}\right)\bigvee\kappa,
\labe{Eq-nrn}
\end{equation}
for some constant $\kappa>1$. Let the observations $X_i$ be i.i.d.\ with a distribution $P$ such that
\begin{equation}
\frac{1}{\nu\left({\scr B}(J_1(n))\right)}\int_{{\scr B}(J_1(n))}K\left(P,P_t\right)\,d\nu(t)=K_n<+\infty.
\labe{Eq-Kullapprox}
\end{equation}
Consider a Bayes estimator $P_{\widetilde{s}}$ of $P$ where $\widetilde{s}$ is a minimizer over $S$ of the function $u\mapsto\Bbb{E}\left[w\left(h(u,\bm{t})\right)\!|\bm{X}\right]$ and that the loss function $w$ satisfies (\ref{Eq-loss}) with $B'\le(\kappa-1)/4$. Then
\[
\Bbb{E}\left[h^2\left(\st\widetilde{s}(\bm{X}),s\right)\right]\le C\Delta^2\left(4^{-J_1(n)}+K_n\right)\quad
\mbox{with }\;\Delta=4\left([(5/4)(1+a'/2)]^{1/\delta}+1\right),
\]
for some universal constant $C$.
\end{theorem}
It immediately follows that 
\begin{equation}
\Bbb{E}\left[h^2\left(P_{\widetilde{s}},P\right)\right]\le 2C\Delta^2\left(4^{-J_1(n)}+K_n\right)
+2h^2(P_s,P).
\labe{Eq-aux8}
\end{equation}
In this bound, $K_n$ and $h^2(P_s,P)$ play the role of bias terms while $4^{-J_1(n)}$ can be viewed as a variance term which depends on both the metric structure of the model $S$ via the function $D$ and the behaviour of $\nu$ in vicinities of $s$. Since $K_n$ depends on $J_1(n)$,  one should choose $J_1(n)$ in order to minimize the risk bound.

An obvious bound for $K_n$ is $K_n\le\kappa_n=\sup_{t\in{\scr B}(J_1(n))}K(P,P_t)$ but, due to the regularization properties of the averaging operation, $K_n$ can be substantially smaller than $\kappa_n$. It can also be controlled, in the case of bounded likelihood ratios, from the following improved version of Lemma~4.4 page 208 of Birg\'e (1983). We include a proof in Section~\ref{P} for the sake of completeness.
%
\begin{lemma}\lab{L-hKull}
For any two probabilities $P$ and $Q$,
\begin{equation}
K(P,Q)\ge-2\log[\rho(P,Q)]\ge2h^2(P,Q).
\labe{Eq-Kull0}
\end{equation}
Moreover, if $\sup_{x\in\scr{X}}(dP/dQ)(x)=M<+\infty$,
\begin{equation}
2\le\frac{K(P,Q)}{h^2(P,Q)}\le {\cal K}(M)\quad
\mbox{with} \;\; {\cal K}(z)=\frac{z(\log z-1)+1}{(z+1)/2-\sqrt{z}}.     
\labe{Eq-Kull1}
\end{equation}
In particular, ${\cal K}(z)\le4+2\log z$ for $z\ge1$.
\end{lemma}
Another consequence of this lemma is the fact that $2h^2(P_s,P)\le K(P,P_s)$. Therefore in typical situations the additional term $2h^2(P_s,P)$ in (\ref{Eq-aux8}) will not be of larger order than $K_n$. 

\subsection{Two preliminary results\labs{C2}}
Let us first consider the situation of a probability $Q$ on some metric space $({\cal S},d)$ of 
diameter $\sup_{t,u\in{\cal S}\times{\cal S}}d(t,u)\le1$ with its Borel $\sigma$-algebra. 
Then the following result, to be proved in Section~\ref{P}, holds:
%
\begin{proposition}\lab{P-Bayesconv}
Let $w$ be some nondecreasing function on $[0,1]$ such that (\ref{Eq-loss}) holds for all $z\in(0,1/2]$ with $\delta,a',B'>0$. Let $J\in\Bbb{N}^*$ and the random variable $\bm{t}$ on $S$ be distributed according to some probability $Q$ which satisfies
\begin{equation}
Q\left[{\cal B}\left(s,2^{j-J}\right)\right]\ge1-a\exp\left[-B4^{j}\right]\quad\mbox{for all }j\in\Bbb{N}
\:\mbox{ with }\:j<J
\labe{Eq-Concen}
\end{equation}
and $a>0$, $B\ge 2$, $ae^{-B}\le1/5$. Then $Q\left[{\cal B}\left(s,2^{-J}\right)\right]\ge4/5$ and 
if  $B'\le (B-1)/4$, any minimizer $\widetilde{u}$ of the function $u\mapsto\Bbb{E}[w(d(\bm{t},u))]$ satisfies
\[
d\left(\widetilde{u},s\right)\le\left(\left[(5/4)(1+0.4aa')\st\right]^{1/\delta}+1\right)2^{-J}.
\]
\end{proposition}

Our next result is a suitable version, due to Yannick Baraud, of a lemma of Barron (1991, Section~5.3) which appears as Proposition~1 in Baraud, Birg\'e and Sart (2014).
%
\begin{proposition}\lab{P-KL}
For any pair $Q,R$ of probabilities such that $K(R,Q)<+\infty$ and any random variable $T$ such that $\Bbb{P}_{Q}[T\ge z]\le ae^{-bz}$ for all $z\ge z_0\ge0$ with $a,b>0$, we have
\[
\Bbb{E}_{R}[T]\le z_0+ b^{-1}\left(1+A+\sqrt{2A}\right)\quad\mbox{with}\quad
A=\log\left(1+ae^{-bz_0}\right)+K(R,Q).
\]
\end{proposition}
%

\subsection{Proof of Theorem~\ref{T-Bayesest}\labs{C4}}
We recall from Section~\ref{A3} that (\ref{Eq-nrn}) and Theorem~\ref{T-Bayes} imply that, for 
$0\le J\le J_1(n)-2$, except on a set $\Gamma_J\subset\scr{X}^n$ with $P_{{\scr B}(J_1(n))}$-probability bounded by $1.05\exp\left[-4^{-J-3}n\right]$, the posterior distribution satisfies (\ref{Eq-aux1}), which can be written
\[
\overline{\nu}\left[\left.{\cal B}\left(s,2^{k-J}\right)\right|\bm{X}\right]\ge1-1.05\exp\left[-4^{k-J-3}n\right]\quad\mbox{for }k\in\Bbb{N}, \;k\le J.
\]
Moreover, $4^{-J-3}n\ge\kappa$ by (\ref{Eq-nrn}). This means that, when $\bm{X}\not\in\Gamma_J$, the posterior distribution $\overline{\nu}(\cdot|\bm{X})$ satisfies the assumptions (\ref{Eq-Concen}) with $a=1.05$ and $B=4^{-J-3}n\ge\kappa$. Therefore, according to Proposition~\ref{P-Bayesconv}, since the function $w$ satisfies (\ref{Eq-loss}) with
\[
B'\le(\kappa-1)/4\le(B-1)/4,
\]
the corresponding Bayes estimator $\widetilde{s}(\bm{X})$ satisfies
\[
h\left(\st\widetilde{s}(\bm{X}),s\right)\le\left([(5/4)(1+1.05\times0.4a')]^{1/\delta}+1\right)2^{-J}
<2^{-J-2}\Delta,
\]
with $\Delta=4\left([(5/4)(1+a'/2)]^{1/\delta}+1\right)>4$. This means that
\[
\Bbb{P}_{{\scr B}(J_1(n))}\left[h\left(\st\widetilde{s}(\bm{X}),s\right)\ge2^{-J-2}\Delta\right]
\le1.05\exp\left[-4^{-J-3}n\right]\quad\mbox{for }0\le J\le J_1(n)-2
\]
or, equivalently,
\[
\Bbb{P}_{{\scr B}(J_1(n))}\left[h^2\left(\st\widetilde{s}(\bm{X}),s\right)\ge4^{-k}\Delta^2\right]
\le1.05\exp\left[-4^{-k-1}n\right]\quad\mbox{for }2\le k\le J_1(n)
\]
and, consequently, since $\Delta>4$,
\begin{equation}
\Bbb{P}_{{\scr B}(J_1(n))}\left[h^2\left(\st\widetilde{s}(\bm{X}),s\right)\ge z\right]\le
1.05\exp\left[-nz/\!\left(16\Delta^2\right)\right]\quad\mbox{for }1\ge z\ge4^{-J_1(n)}\Delta^2.
\labe{Eq-Baycon}
\end{equation}
This allows us to apply Proposition~\ref{P-KL} with $Q=P_{{\scr B}(J_1(n))}$, $R=
P^{\otimes n}$, $T=h^2\left(\st\widetilde{s}(\bm{X}),s\right)$, $z_0=4^{-J_1(n)}\Delta^2$, $a=1.05$ and $b=n/\!\left(16\Delta^2\right)$, from which we derive that
\begin{equation}
\Bbb{E}\left[h^2\left(\st\widetilde{s}(\bm{X}),s\right)\right]\le4^{-J_1(n)}\Delta^2+ 
16\Delta^2n^{-1}\left(1+A+\sqrt{2A}\right)
\labe{Eq-Bay}
\end{equation}
with
\begin{eqnarray*}
A&=&\log\left(1+1.05\exp\left[-4^{-J_1(n)-2}n\right]\right)+K\left(P^{\otimes n},
P_{{\scr B}(J_1(n))}\right)\\&<&1.05\exp\left[-4^{-J_1(n)-2}n\right]+K\left(P^{\otimes n},
P_{{\scr B}(J_1(n))}\right).
\end{eqnarray*}
Since $K(P,Q)=\int\log(dP/dQ)dP$ and $-\log$ is a convex function, 
\begin{eqnarray*}
K\left(P^{\otimes n},P_{{\scr B}(J_1(n))}\right)&\le&\frac{1}{\nu\left({\scr B}(J_1(n))\right)}
\int_{{\scr B}(J_1(n))}K\left(P^{\otimes n},P_t^n\right)\,d\nu(t)\\&=&
\frac{n}{\nu\left({\scr B}(J_1(n))\right)}\int_{{\scr B}(J_1(n))}K\left(P,P_t\right)\,d\nu(t)\;\;=\;\;nK_n
\end{eqnarray*}
and our conclusion follows.

\Section{Proofs\labs{P}}

\subsection{Proof of Lemma~\ref{L-bayconv}\labs{P1}}
We may assume, without loss of generality, that $\mu$ dominates $P$ with $dP/d\mu=g$. Since $\lambda$ is a probability, Jensen's Inequality implies that $\sqrt{f_\lambda(x)}\ge
\int_S\sqrt{f_t(x)}\,d\lambda(t)$. If $h(P,P_t)\le r$, then $\rho(P,P_t)\ge1-r^2$ and, according to Fubini's Theorem,
\begin{eqnarray*}
\rho(P,P_\lambda)&=&\int_{\scr X}\sqrt{g(x)f_\lambda(x)}\,d\mu(x)\;\;\ge\;\;\int_{\scr X}
\sqrt{g(x)}\left(\int_S\sqrt{f_t(x)}\,d\lambda(t)\right)\,d\mu(x)\\&=&
\int_S d\lambda(t)\int_{\scr X}\sqrt{g(x)f_t(x)}\,d\mu(x)\;\;=\;\;
\int_S\rho(P,P_t)\,d\lambda(t)\;\;\ge\;\;1-r^2.\cqfd
\end{eqnarray*}

\subsection{Likelihood ratio tests\labs{P2}}
Our results are based on the following fundamental result which is due to Le Cam and 
can be found in his book, Le Cam (1986, Lemma~2 page~477).
%
\begin{proposition}\lab{P-bayconv}
Let $\{P_t,t\in S\}$ be a dominated  family of probability measures on ${\scr X}$ with densities $f_t$ with respect to some mesure $\mu$. We consider two distributions $Q_1$, $Q_2$ on ${\scr X}$, two positive numbers $r_1,r_2$ such that $h(Q_1,Q_2)>r_1+r_2$ and two probability distributions $\lambda_1,\lambda_2$ on $S$ such that $\lambda_j
\left(\{t\,|\,h(P_t,Q_j)\le r_j\}\right)=1$ for $j=1,2$. We define the distributions $R_1,R_2$ on ${\scr X}^n$ by their densities with respect to $\mu^{\otimes n}$, 
\[
\frac{dR_j}{d\mu^{\otimes n}}(\bm{x})=\int_Sf_{t}(\bm{x})\,d\lambda_j(t)\quad\mbox{for all } \bm{x}\in{\scr X}^n\;\;\mbox{and}\;\;j=1,2\quad\mbox{with }
f_t(\bm{x})=\prod_{i=1}^nf_t(x_i).
\]
Then 
\[
\rho(R_1,R_2)\le\left[1-\left(\st h(Q_1,Q_2)-r_1-r_2\right)^2\right]^n.
\]
\end{proposition}
%
\begin{corollary}\lab{C-Main}
Let $\{P_t,t\in S\}$ be a dominated  family of probability measures on ${\scr X}$ with densities $f_t$ with respect to some mesure $\mu$ and let $\nu$ be a probability on $S$. Let $B_1$ and $B_2$ be two measurable subsets of $S$ such that $\nu(B_1)\nu(B_2)>0$ and 
$B_j\subset{\cal B}_j$ for $j=1,2$ where ${\cal B}_1$ and ${\cal B}_2$ are two Hellinger balls in $S$ such that $h({\cal B}_1,{\cal B}_2)>0$. Then, given a random element $\bm{X}=(X_1,\ldots,X_n)\in{\scr X}^n$ with density $g_{B_1}$ with respect to $\mu^{\otimes n}$ given by (\ref{Eq-gB}), 
\begin{eqnarray}
\Bbb{P}_{B_1}\left[\log\left(\frac{g_{B_2}(\bm{X})}
{g_{B_1}(\bm{X})}\right)\ge y\right]&\le&e^{-y/2}\rho\left(P_{B_1},P_{B_2}\right)\nonumber\\
&\le&\exp\left[-\frac{y}{2}-nh^2({\cal B}_1,{\cal B}_2)\right]\quad\mbox{for all }y\in\Bbb{R}.
\labe{Eq-Lap2}
\end{eqnarray}
\end{corollary}
%
%
{\em Proof:}
The first inequality is an application of Lemma~\ref{L-expineq} and the second follows from Proposition~\ref{P-bayconv} applied with $\lambda_j=P_{B_j}$ and ${\cal B}_j=
\{t\,|\,h(P_t,Q_j)\le r_j\}$ so that $h({\cal B}_1,{\cal B}_2)\ge h(Q_1,Q_2)-r_1-r_2$.\cqfd
%
\subsection{Proof of Theorem~\ref{T-Bayes}\labs{P3}}
Let us first observe that the monotonicity of $D$ and $\beta$ and the fact that $\gamma\le4$
imply that condition (\ref{Eq-Required}) can only be satisfied for some pair $(J_1,J)$ if it is satisfied for $J_1=3$ and $J=1$ which shows that (\ref{Eq-n}) is a necessary and sufficient condition for the existence of such a pair $(J_1,J)$.

Given $J_1$ and $J$ that satisfy (\ref{Eq-Required}), we can partition $\left[{\scr B}(J)\right]^c$ into sets $F_{l,k}$, $0\le l\le J-1$, $1\le k\le K_l$ in the following way: we set $F_l=
{\scr B}(l)\setminus{\scr B}(l+1)$ and consider a maximal $\left(2^{-l-2}\right)$-separated subset $N_l$ of $F_l$, hence of ${\scr B}(l)$. By Assumption~\ref{A-Dim}, 
\begin{equation}
K_l=|N_l|\le\exp\left[D\!\left(2^{-l-2}\right)\right].
\labe{Eq-Kj}
\end{equation}
This induces a covering of $F_l$ by $K_l$ balls ${\cal B}_{l,k}$, $1\le k\le K_l$ with radius $2^{-l-2}$ and centers in $N_l$ from which one can deduce a partition of $F_l$ into $K_l$ subsets $F_{l,k}$, each one being included in the corresponding ball ${\cal B}_{l,k}$ with center in $F_l$. It follows that
\[
h\!\left({\scr B}(J_1),{\cal B}_{l,k}\right)\ge h(s,F_l)-2^{-l-2}-2^{-J_1}=2^{-l-1}-2^{-l-2}-2^{-J_1}
\;\;\mbox{ for }1\le k\le K_l
\]
and, since $l\le J-1\le J_1-3$ by (\ref{Eq-Required}),
\[
h\!\left({\scr B}(J_1),{\cal B}_{l,k}\right)\ge2^{-l-2}-2^{-J_1}=
\left(1-2^{-J_1+l+2}\right)2^{-l-2}\ge2^{-l-3}.
\]
Then, since $F_{l,k}\subset{\cal B}_{l,k}$, it follows from (\ref{Eq-Lap2}) that
\[
\Bbb{P}_{{\scr B}(J_1)}\left[\log\left(\frac{g_{F_{l,k}}(\bm{X})}{g_{{\scr B}(J_1)}(\bm{X})}\right)\ge y\right]
\le\exp\left[-\frac{y}{2}-4^{-l-3}n\right]\quad\mbox{for all }y\in\Bbb{R}.
\]
Setting $y=-4^{-l-3}n$, we get
\[
\Bbb{P}_{{\scr B}(J_1)}\left[\log\left(\frac{g_{F_{l,k}}(\bm{X})}{g_{{\scr B}(J_1)}(\bm{X})}\right)
\ge-4^{-l-3}n\right]\le\exp\left[-(2n)\times4^{-l-4}\right].
\]
This means that
\begin{equation}
\frac{g_{F_{l,k}}(\bm{X})}{g_{{\scr B}(J_1)}(\bm{X})}<\exp\left[-4^{-l-3}n\right]
\quad\mbox{for }0\le l\le J-1\mbox{ and }1\le k\le K_l,
\labe{Eq-A3}
\end{equation}
except if $\bm{X}\in\Gamma_J$ with
\[
P_{{\scr B}(J_1)}[\Gamma_J]\le\sum_{l=0}^{J-1}K_l\exp\left[-(2n)\times4^{-l-4}\right]\le
\sum_{l=0}^{J-1}\exp\left[-(2n)\times4^{-l-4}+D\!\left(2^{-l-2}\right)\right],
\]
by (\ref{Eq-Kj}). Since $1\le l+1\le J$, (\ref{Eq-Req}) applies with $J'=l+1$ leading to
\begin{equation}
4^{-l-4}n\ge\left[\gamma^{J_1-l}\beta(J_1)/3\right]\bigvee D\!\left(2^{-l-2}\right)\ge1
\quad\mbox{for }0\le l\le J-1.
\labe{Eq-Re}
\end{equation}
In particular, $D\!\left(2^{-l-2}\right)\le4^{-l-4}n$, which shows that $P_{{\scr B}(J_1)}[\Gamma_J]\le\sum_{l=0}^{J-1}\exp\left[-4^{-l-4}n\right]$. Then (\ref{Eq-T1a}) follows since
\begin{eqnarray}
\lefteqn{\sum_{l=0}^{J-1}\exp\left[-4^{-l-4}n\right]}\hspace{10mm}\nonumber\\&<&\sum_{k\ge0}
\exp\left[-4^{-J-3+k}n\right]\;\;\le\;\;\exp\left[-4^{-J-3}n\right]\sum_{k\ge0}
\exp\left[-4^{-J-3}n\left(4^k-1\right)\right]\nonumber\\&\le&\exp\left[-4^{-J-3}n\right]
\sum_{k\ge0}\exp\left[-\left(4^k-1\right)\right]\;\;<\;\;1.05\exp\left[-4^{-J-3}n\right],
\labe{Eq-series}
\end{eqnarray}
where we used the fact that $4^{-J-3}n\ge1$.

Let us now turn to the proof of (\ref{Eq-T1b}). The posterior distribution $\overline{\nu}(A|\bm{X})$ of any subset $A\supset{\scr B}(J_1)$ of $S$ is given by
\begin{eqnarray*}
\frac{1}{\overline{\nu}(A|\bm{X})}&=&\int_Sf_t(\bm{X})\,d\nu(t)
\left[\int_{A}f_t(\bm{X})\,d\nu(t)\right]^{-1}\\&=&
\left[\int_{A}f_t(\bm{X})\,d\nu(t)+\int_{A^c}f_t(\bm{X})\,d\nu(t)\right]
\left[\int_{A}f_t(\bm{X})\,d\nu(t)\right]^{-1}\\&\le&1+\int_{A^c}f_t(\bm{X})\,d\nu(t)
\left[\int_{{\scr B}(J_1)}f_t(\bm{X})\,d\nu(t)\right]^{-1},
\end{eqnarray*}
so that
\begin{equation}
\overline{\nu}(A|\bm{X})\ge1-
\left[\nu\!\left(\st{\scr B}(J_1)\right)g_{{\scr B}(J_1)}(\bm{X})\right]^{-1}\int_{A^c}f_t(\bm{X})\,d\nu(t).
\labe{Eq-A1}
\end{equation}
For $1\le j\le J$ (since there is nothing to prove for $j=0$),
\begin{equation}
\int_{[{\scr B}(j)]^c}f_t(\bm{X})\,d\nu(t)=\sum_{l=0}^{j-1}\,\sum_{k=1}^{K_l}
\int_{F_{l,k}}f_t(\bm{X})\,d\nu(t)=\sum_{l=0}^{j-1}\,\sum_{k=1}^{K_l}\nu(F_{l,k})g_{F_{l,k}}(\bm{X}).
\labe{Eq-A2}
\end{equation}
Putting (\ref{Eq-A1}), (\ref{Eq-A2}) and (\ref{Eq-A3}) together, we derive that, when 
$\bm{X}\not\in\Gamma_J$,
\begin{eqnarray*}
\overline{\nu}\left[\left.\left({\scr B}(j)\st\right)^c\right|\bm{X}\right]&\le&
\frac{\sum_{l=0}^{j-1}\sum_{k=1}^{K_l}\nu(F_{l,k})
g_{F_{l,k}}(\bm{X})}{\nu\!\left(\st{\scr B}(J_1)\right)g_{{\scr B}(J_1)}(\bm{X})}\\&<&\frac{1}
{\nu\!\left(\st{\scr B}(J_1)\right)}\sum_{l=0}^{j-1}\,\sum_{k=1}^{K_l}\nu(F_{l,k})
\exp\left[-4^{-l-3}n\right]\\&=&\sum_{l=0}^{j-1}\frac{\nu(F_{l})}{\nu\!\left(\st{\scr B}(J_1)\right)}
\exp\left[-4^{-l-3}n\right]\;\;\le\;\;\sum_{l=0}^{j-1}\frac{\nu\!\left({\scr B}(l)\right)}
{\nu\!\left(\st{\scr B}(J_1)\right)}\exp\left[-4^{-l-3}n\right]\\&\le&\sum_{l=0}^{j-1}
\exp\left[\gamma^{(J_1-l)}\beta(J_1)-4^{-l-3}n\right],
\end{eqnarray*}
where the last inequality follows from Assumption~\ref{A-Prior}. Applying (\ref{Eq-Re}) again, we derive that $\overline{\nu}\left[\left.\left({\scr B}(j)\st\right)^c\right|\bm{X}\right]\le
\sum_{l=0}^{j-1}\exp\left[-4^{-l-4}n\right]$ and conclude as we did for (\ref{Eq-series}).

\subsection{Proof of Proposition~\ref{P-Bayesconv}\labs{P4}}
Let us first evaluate $\Bbb{E}\left[w(d(\bm{t},s))\right]$. Since $B-4B'\ge1$,
\begin{eqnarray*}
\Bbb{E}[w(d(\bm{t},s))]
&\le&w\left(2^{-J}\right)+\sum_{j=0}^{J-1}w\left(2^{j+1-J}\right)\Bbb{P}\left[2^{j-J}<d(\bm{t},s))\le
2^{j+1-J}\right]\\&\le&w\left(2^{-J}\right)+a\sum_{j=0}^{J-1}w\left(2^{j+1-J}\right)\exp\left[-B4^{j}\right]\\&\le&
w\left(2^{-J}\right)\left[1+aa'\sum_{j=0}^{J-1}\exp\left[-B4^{j}+B'4^{j+1}\right]\right]\\&\le&
w\left(2^{-J}\right)\left[1+aa'\sum_{j=0}^{J-1}\exp\left[-4^{j}\right]\right]\;\;\le\;\;(1+0.4aa')w\left(2^{-J}\right).
\end{eqnarray*}
Let now $t_1\in{\cal S}$ with $d(t_1,s)\ge[2^{j}+1]2^{-J}$. If $d(\bm{t},s)\le2^{-J}$, then 
$d(\bm{t},t_1)\ge2^{j-J}$ hence,
\[
\Bbb{E}\left[w(d(\bm{t},t_1))\right]\ge w\left(2^{j-J}\right)
Q\left[d(\bm{t},s)\le2^{-J}\right]\ge(4/5)2^{j\delta}w\left(2^{-J}\right).
\]
It follows that $t_1$ cannot be a minimizer of the function $t\mapsto 
\Bbb{E}\left[w(d(\bm{t},t))\right]$ as soon as 
\[
(4/5)2^{j\delta}>(1+0.4aa')\quad\mbox{ or }\quad 2^j>[(5/4)(1+0.4aa')]^{1/\delta}
\]
and the conclusion follows.

\subsection{Proof of Lemma~\ref{L-hKull}\labs{P5}}
To prove (\ref{Eq-Kull0}) we use Jensen's Inequality,
\[
-\frac{1}{2}K(P,Q)=\int\log\left(\sqrt{\frac{dQ}{dP}}\right)dP\le
\log\left(\int\sqrt{dQdP}\right)=\log[\rho(P,Q)],
\]
and $\log\left(1-h^2(P,Q)\right)\le-h^2(P,Q)$. 
Let $\eta=dP/dQ$ so that $\int\eta\,dQ=1$. Then
\begin{eqnarray*}
K(P,Q)&=&\int\log\eta\,dP\;\;=\;\;\int\eta\log\eta\,dQ\\&=&
\int(\eta\log\eta-\eta+1)\,dQ\;\;=\;\;\int{\cal K}(\eta)
\left[(\eta+1)/2-\sqrt{\eta}\right]dQ,
\end{eqnarray*}
while
\[
h^2(P,Q)=1-\int\sqrt{\eta}\,dQ=\int\left[(\eta+1)/2-\sqrt{\eta}\right]dQ.
\]
Since the function ${\cal K}$ is continuous and increasing on $\Bbb{R}_+$ with ${\cal K}(0)=2$, it
follows that $2\le{\cal K}(\eta)\le{\cal K}(M)$ hence (\ref{Eq-Kull1}). The bound on ${\cal K}$
follows from calculus.\vspace{3mm}\\
{\large{\bf Acknowledgements}}
I would like to thank Judith Rousseau for her invitation to participate to the Bayesian Nonparametrics Conference 2013, which strongly pushed me to write down properly some old sketches about Bayesian procedures and also Yannick Baraud, whose numerous criticisms and suggestions highly contributed to the improvement of the presentation of the results.
\vspace{6mm}\\
%
%
\Large{\bf {References\vspace{1mm}}}
\setlength{\parskip}{1mm}
{\small

BARAUD, Y., BIRG\'E, L.  and SART, M. (2014). A new method for estimation and model selection: $\rho$-estimation. ArXiv e-print 1403.6057.

BARRON, A.R. (1991). Complexity regularization with applications to artificial neural 
networks. In {\it  Nonparametric Functional Estimation} (G. Roussas, ed.). Kluwer,
Dordrecht, 561-576. 

BIRG\'E, L.  (2006). Model selection via testing : an alternative to (penalized) maximum
likelihood estimators. {\it Ann. Inst. Henri Poincar\'e, Probab.\ et Statist.} {\bf 42}, 273-325.

GHOSAL, S., GHOSH, J.K. and van der VAART, A.W. (2000). Convergence rates of posterior
distributions. {\it Ann. Statist.}  {\bf28}, 500-531.  

GHOSAL, S., LEMBER and van der VAART, A.W. (2003). On Bayesian adaptation. {\it Acta
Applicandae Mathematicae}  {\bf79}, 165-175.  

IBRAGIMOV, I.A. and HAS'MINSKII, R.Z.  (1981). {\it Statistical Estimation: Asymptotic 
Theory}. Springer-Verlag, New York.

KOLMOGOROV, A.N. and TIKHOMIROV, V.M. (1961). $\varepsilon $-entropy and   
$\varepsilon$-capacity of sets in function spaces. {\it  Amer. Math. Soc. Transl. (2)} {\bf17},
277-364.  

Le CAM, L.M. (1973). Convergence of estimates under dimensionality restrictions. {\it
Ann. Statist.}  {\bf1} , 38-53. 

Le CAM, L.M. (1982). On the risk of Bayes estimates. {\it Statistical Decision Theory and
Related Topics III, Vol.\ 2} (S.S.~Gupta and J.O.~Berger, eds.), 121-137.  Academic Press, 
New York. 

Le CAM, L.M. (1986). {\it Asymptotic Methods in Statistical Decision Theory}. 
Springer-Verlag, New York. 

LORENTZ, G.G. (1966). {\it Approximation of Functions}. Holt, Rinehart, Winston,  New
York. 

MASSART, P.  (2007). Concentration Inequalities and Model Selection. In {\it  Lecture on
Probability Theory and Statistics, Ecole d'Et\'e de  Probabilit\'es de Saint-Flour XXXIII -
2003} (J.~Picard, ed.). Lecture Note in Mathematics, Springer-Verlag,  Berlin. 

van der VAART, A.W. (2003).  {\it Convergence rates of nonparametric posteriors. Highly Structured Stochastic Systems} P.J. Green, N. Hjort and S. Richardson (eds). Oxford University Press, Oxford.
\vspace{8mm}\\
Lucien BIRG\'E\\ UMR 7599 ``Probabilit\'es et mod\`eles al\'eatoires"\\
Laboratoire de Probabilit\'es, bo\^{\i}te 188\\
Universit\'e Paris 06, 4 Place Jussieu\\
F-75252 Paris Cedex 05\\
France\vspace{2mm}\\
e-mail: lucien.birge@upmc.fr

\mbox{}%
}

\end{document}